\documentclass[10pt,leqno]{article}

\usepackage{amsmath,amssymb,amsthm,mathrsfs,dsfont}

\usepackage[margin=3cm]{geometry} 

\usepackage{titlesec,hyperref}

\usepackage{color}

\usepackage{fancyhdr}
\titleformat{\subsection}{\it}{\thesubsection.\enspace}{1pt}{}

\newtheorem{theo}{Theorem}[section]

\newtheorem{rema}[theo]{Remark}
\numberwithin{equation}{section}

\allowdisplaybreaks 

\begin{document}
\title{Generic singularity behavior of conservative solutions to the Novikov equation
\hspace{-4mm}
}

\author{ Zhen $\mbox{He}^1$ \footnote{E-mail:hezh56@mail2.sysu.edu.cn},\quad
	Wei $\mbox{Luo}^1$\footnote{E-mail:luowei23@mail2.sysu.edu.cn  } \quad and\quad
	Zhaoyang $\mbox{Yin}^{1,2}$\footnote{E-mail:mcsyzy@mail.sysu.edu.cn}\\
	$^1\mbox{Department}$ of Mathematics,
	Sun Yat-sen University, Guangzhou 510275, China\\
	$^2\mbox{School}$ of Science,\\ Shenzhen Campus of Sun Yat-sen University, Shenzhen 518107, China}

\date{}
\maketitle
\hrule

\begin{abstract}
In this paper, we concentrate on the Novikov equation. We provide a description of the solution in a neighborhood of each
singular point. The result shows the difference of the conservative solution's singularity behavior the  between the Camassa-Holm equaation and the Novikov equation  \\
\vspace*{5pt}
\noindent {\it 2010 Mathematics Subject Classification}: 35Q30, 76B03, 76D05, 76D99.

\vspace*{5pt}
\noindent{\it Keywords}: Novikov equation; generic regularity; conservative weak solution;singularity
\end{abstract}

\vspace*{10pt}

\section{Introduction}

   Consideration here is the initial-value problem for the Novikov equation in the form
 \begin{equation}\label{001}
 	\left\{\begin{array}{l}
 		u_t+u^2u_x+\partial_xP_1+P_2=0, \\
 		u(0,x)=u_0(x),  \quad x\in\mathbb{R}.
 	\end{array}\right.
 \end{equation}
 where
 $$P_1\triangleq p\ast(\frac{3}{2}u u_x^2+u^3)\quad \textit{and} \quad P_2\triangleq\frac{1}{2}p\ast u_x^3$$
 with $p(x)=\frac{1}{2}e^{-|x|}$.
  The equation was proposed by Novikov in \cite{NE}. In \cite{Chen-Liu}, Chen, Hu and Liu showed that the Novikov equation can be viewed as a shallow water model.

    There is a huge literature devoted to the equations \eqref{001}, local well-posedness can be referred to
 \cite{GY,HH,LY,WY,WY2,YL,YL2}.  Chen, Zhang and Liu\cite{CWN} proved the existence and uniqueness of conservative solutions for the Novikov equation in $H^1\cap W^{1,4}$.

  Recently, Li and Zhang \cite{GRCH} proved the generic property and the singular behavior of the Camassa-Holm equation and the two-component Camassa-Holm equation. The original result was studied by Schaeffer in \cite{GRCL}, which showed that for the one space dimensional conservation law, the generic solutions are piecewise smooth, with finitely shocks in a bounded domain in the $(t,x)$ plane. Bressan and his collaborators\cite{GRNW,GRDW,CSNW} studied generic property and singularity behavior for the variational wave equation. Yang \cite{Yang} studied generic regularity of energy conservative solutions to the rotation Camassa-Holm equation. Cai, Chen, Shen and Tan \cite{Cai} studied the generic property of conservative solutions to the Hunter-Saxton type equations and give a new way to construct a Finsler type metric which renders the flow
  uniformly Lipschitz continuous on bounded subsets of $H^1(\mathbb{R^+})$.

  The generic property of the Novikov equation has  been studied in \cite{NOV}. In \cite{CWN}, the authors pointed out that the singularity for the conservative solutions of the Novikov equation is difference between the Camassa-Holm equation. The wave speed $c(u)$ for the Camassa-Holm equation is $u$ and for the Novikov equation is $u^2$. Because of the difference nonlinearity of wave speed, the singularity behaviour is various. In this paper, we give a exactly proof to show that why the singularity is various between the Camassa-Holm equation and the Novikov equation.

  From the above generic regularity, we can obtain the asymptotic description of the solution in a neighborhood of each singular point, where $|u_x|\rightarrow\infty$.
  \begin{theo}
  	Consider generic initial data $u_0\in C^3(\mathbb{R}) \cap H^1(\mathbb{R})
  	\cap W^{1,4}(\mathbb{R})$ as in \cite{NOV}, with $u_0\in C^\infty(\mathbb{R})$. Call $(u,v,\xi,x,t)$ the corresponding solution of the semilinear system \eqref{006} and let u=u(x,t) be the solution to the original equation \eqref{001}. Consider a singular point $P=(t_0,Y_0)$ where v=$\pi$, and set $(x_0,t_0)=(x(t_0,Y_0),t(t_0,Y_0))$. Generically, at the singular point, u has following parameteric expression.
  	\begin{itemize}
  		\item [\textcircled{\scriptsize{1}}]If P is a point of Type $\mathcal{I}$ , i.e.$v_Y=0$ and $v_{YY}\neq 0$ $v=\pi,~~v_Y$ then
  		\begin{equation}\label{S1}
  			u(t,x)=A(x-x_0)^{\frac{3}{4}}+B(t-t_0)+\mathcal{O}(1)(| t-t_0 |^2+| x-x_0|^{\frac{7}{8}})
  		\end{equation}
  		for some constant $A, B$.
  		\item [\textcircled{\scriptsize{2}}]If P is a point of Type $\mathcal{I}\mathcal{I}$ , i.e.$v_Y\neq 0$ and $v_{YY}= 0$ $v=\pi,~~v_Y$ then
  		\begin{equation}\label{S2}
  			u(t,x)=A(x-x_0)^{\frac{4}{5}}+B(t-t_0)+\mathcal{O}(1)(| t-t_0 |^2+| x-x_0|)
  		\end{equation}
  		for some constant $A, B$.
  	\end{itemize}
  	
  \end{theo}

\section{Generic singularity behavior}
~~~For smooth data $u_0\in C^\infty(\mathbb{R})$, the solution $(t,Y)\rightarrow(x,t,u,v,q)(t,Y)$ of the semilinear system \eqref{006}, with initial data as in \eqref{i2}, remains smooth on the entire t-Y plane. Yet the smoothness of the solution u of \eqref{001} is still needed to study because the coordinate change:(Y,t)$\rightarrow$(x,t) is not smoothly invertible. By definitions, its Jacobian matrix is computed by
\begin{align}
	\begin{pmatrix}
		x_Y & x_t
		\\
		t_Y & t_t
	\end{pmatrix}
=\begin{pmatrix}
	qcos^4\frac{v}{2} & u
	\\
	0 & 1
\end{pmatrix}.
\end{align}
  And we will observe that the matrix is invertible when $v\neq\pi$. To study the set of points in the $t-x$ plane where $u$ is singular, we thus need to look at points where $v=\pi$.

   For simplicity, we shall assume that the initial data $u_0$ are smooth , so we shall not need to count how many derivatives are actually used to derive the Talyor approximations.

\begin{proof}
	
	Following the idea in \cite{CWN}, the characteristic equation is
	\begin{align}
		\frac{dx(t)}{dt}=u^2(t,x(t)).
	\end{align}
	
	If we consider $Y=Y(t,x)$ is a characteristic coordinate, and denote $T=t$.
	Then we consider function  $f(t,x)=f(T,x(t,Y)) $ as a function of $(T,Y)$ also denoted by $ f(T,Y)$.It is easy to check that
	\begin{align}
		f_t+u^2f_x=f_Y(Y_t+u^2Y_x)+f_t(T_t+u^2T_x)=f_T.
	\end{align}
	We also denote
	\begin{align}\label{zh1}
		v=2\arctan u_x , ~~~q=\frac{(1+u_x^2)^2}{Y_x}.                                                                                                                                                                                                                   \end{align}
	
	Then the  conservative solution is constructed by the following semilinear system
	\begin{equation}\label{006}
		\left\{\begin{array}{l}
			u_T=-\partial_xP_1-P_2  \\
			v_T=-u\sin^2\frac{v}{2}+2u^3\cos^2\frac{v}{2}-2\cos^2\frac{v}{2}(P_1+\partial_xP_2),  \\
			q_T=q[(2u^3+u)-2(P_1+\partial_xP_2)]\sin v,
		\end{array}\right.
	\end{equation}
	with the initial condition
	\begin{equation}\label{i2}
		\left\{\begin{array}{l}
			u(0,\beta)=u_0(x(0,\beta))  \\
			v(0,\beta)=2\arctan(u_0^{\prime}(x(0,\beta))),  \\
			q(0,\beta)=1.
		\end{array}\right.
	\end{equation}
	for every $\beta \in \mathbb{R}$.
	
	To study the singularities of the solution $u$ of \eqref{001}, we should focus on the level sets $\{v(t,Y)=\pi\}.$

	Let $(u,v,\xi)$ be the smooth solution of \eqref{006}. And take derivatives to the equation of $v$, we obtain
	\begin{align}
		\frac{\partial}{\partial_T}v_Y=&-u_Y\sin^2(\frac{v}{2})-u\cos(\frac{v}{2})\sin(\frac{v}{2})v_Y+6u_Yu^2\cos^2(\frac{v}{2})-2u^3(\sin\frac{v}{2}\cos\frac{v}{2}v_Y)
		\notag\\
		&+\cos(\frac{v}{2})\sin(\frac{v}{2})v_Y(P_1+\partial_xP_2)-2\cos^2\frac{v}{2}(\partial_YP_1+\partial_Y\partial_xP_2)
		\notag \\
		&=-u_Y\sin^2(\frac{v}{2})-\frac{1}{2}uv_Y\sin v+6u_Yu^2\cos^2(\frac{v}{2})-u^3v_Y\sin v
		\notag \\
		&+\frac{1}{2}\cos vv_Y(P_1+\partial_xP_2)-(\cos v+1)(\partial_YP_1+\partial_Y\partial_xP_2).
	\end{align}
	Following the same line
	\begin{align}
		\begin{split}
			\frac{\partial}{\partial_T}v_{YY}=
			&u_{YY}(\sin^2\frac{Y}{2}+\frac{1}{2}v_Y\sin v+6u^2\cos^2\frac{v}{2}-u^3\sin v)+v_Y(-\frac{1}{2}u_Y\sin v-\frac{1}{2}v_Y\cos v-u^3v_Y\cos v)
			\\
			&-\frac{1}{2}\sin vv_Y^2(P_1+\partial_xP_2)+v_Y(\partial_YP_1+\partial_Y\partial_xP_2)(\frac{1}{2}\cos v+\sin v)-(\cos v+1)(\partial_Y\partial_xP_1+\partial_Y^2\partial_xP_2).
			\\
			&	\frac{\partial}{\partial_t}q_Y= \xi_Y[(2u^3+u)-2(P_1+\partial_xP_2)]\sin v
			+q[(2u^3+u)-2(P_1+\partial_xP_2)]\cos vv_Y
			\\
			&~~~~~~~~~~~~~~~~~~~~~~~~+q[(6u^2u_Y+u_Y)-2(\partial_YP_1+\partial_Y\partial_xP_2)]\sin v,
		\end{split}
	\end{align}
   We consider the situations as 
	1. Let P be the point of Type $\mathcal{I}$  Recalling \eqref{zh1} and \eqref{006} then we will have
	\begin{align}
		u_Y=u_xx_Y=\frac{1}{2}q \cdot \sin~v\cdot \cos^2\frac{v}{2}.
	\end{align}
    In a similar way, we obtain
    \begin{align}
    	u_{YY}&=\frac{1}{2}q_Y\cdot \sin v\cdot \ cos^2\frac{v}{2}+\frac{1}{2}q \cdot v_Y \cdot \cos v \cdot \cos^2\frac{v}{2}-\frac{v_Y}{4}\cdot q \sin^2v,
    	\notag \\
    	u_{Yt}&=\frac{1}{2}q sin^2vcos^2\frac{v}{2}\big((2u^3+u)sin^2\frac{v}{2}-2(P_1+\partial_xP_2)+(2u^3+u)cos^2\frac{v}{2}\big)
    	\notag \\
    	&+\frac{1}{2}(1-4sin^2\frac{v}{2})cos^2\frac{v}{2}\big(-usin^2\frac{v}{2}+2u^3cos^2\frac{v}{2}-2cos^2\frac{v}{2}(P_1+\partial_xP_2)\big).
    \end{align}
   At the point P, provides us with that $v_Y=0, ~v_t=0$ and $v_{YY}\neq 0$.
   If we denote $u_{Y^n}$ that u is differentiated with Y by n times. It is not hard to check
   \begin{align}
   	u_{Y^3}=0,~~~u_{Y^4}=0,~~~u_{Y^5}=0,~~~u_{Y^6}=\frac{6}{8}q v_{YY}^3cosv sin^2\frac{v}{2}-3v_{YY}^3q cos^2v
   \end{align}
   which means $u_{Y^6}\neq0$
   Following the same method, we obtain that
   \begin{align}
   	u_{Y^2t}=0,~~~u_{Y^3t}=0,~~~u_{Y^4t}=0,~~~u_{Y^5t}=\frac{6}{8}q v_{YY}^2v_{Yt}cosv sin^2\frac{v}{2}-3v_{YY}^2v_{Yt}q cos^2v
   \end{align}
   we will have $u_{Y^5t}=0$.

So we have Taylor approximations of u at the singular point $P=(t_0,Y_0)$
\begin{equation}\label{u1}
	u(t,Y)=B_1(t-t_0)+B_3(Y-Y_0)^6+\mathcal{O}(1)( {\lvert
	 t-t_0
	 \rvert} ^2
	  , {\lvert Y-Y_0 \rvert} ^7 )
 \end{equation}
By \eqref{zh1} and \eqref{006}, we can obtain
\begin{align}
	x_Y=\frac{\xi}{(1+u_x^2)^2}=q cos^4\frac{v}{2},
\end{align}
similarly we have
\begin{align}
	x_{YY}&= q_Ycos^4\frac{v}{2}-2q v_Ycos^3\frac{v}{2}sin\frac{v}{2}~~~~~~	x_{Yt}=q_tcos^4\frac{v}{2}-2q, v_tsin\frac{v}{2}cos^3\frac{v}{2},
\end{align}
so it is easy to check that in the point P
$$x_{Y^i}=0,~~~~(i=1,2,3,4,5,6,7)$$
and
\begin{align}
	x_{Y^8}&=36q v_{YY}^4sin^4\frac{v}{2}
~~~~~~
	x_{Y^7t}=36q v_{YY}^3v_{Yt}sin^4\frac{v}{2}.
\end{align}
So we have the Taylor approximations of x at the singular point $P=(t_0,Y_0)$
\begin{equation}\label{x1}
   x(t,Y)=x(t_0,Y)+A_2(Y-Y_0)^8+\mathcal{O}(1)( \lvert t-t_0\rvert^2,\lvert Y-Y_0\rvert^9 ).
\end{equation}
We combine \ref{u1} and \eqref{x1} to deduce \eqref{S1}.

2.If P is of type $\mathcal{I}\mathcal{I}$,we have
\begin{align}
	v=\pi, ~~~v_Y\neq 0,~~~ v_{YY}=0
\end{align}
which implies
\begin{align}
	&u_Y=u_{YY}=u_{Y^3}=u_{Yt}=u_{Y^2t}=u_{Y^3t}=0
	\\
	&u_{Y^4}=\frac{1}{4} q v_Y^3cosv \sin^2\frac{v}{2}-\frac{{v_Y}^3}{2}q\cos^2v=\frac{1}{2}v_Y^3q\neq 0.
\end{align}
Then we obtain the Taylor approximation of u is
\begin{align}\label{u2}
	u(t,Y)=B_1(t-t_0)+B_2(Y-Y_0)^4+\mathcal{O}(1)(|t-t_0|^2+|Y-Y_0|^5)
\end{align}
 At point P,
 \begin{align}
 	&x_Y=0,~~~x_{Y^2}=0,~~~x_{Y^3}=0,~~~x_{Y^4}=0,
 	\\
 	&x_{Yt}=0,~~~x_{Y^2t}=0,~~~x_{Y^3t}=0~~~,~~~x_{Y^4t}=0,
 	\\
 	&x_{Y^5}=\frac{3}{2}q v_Y^3cos^3v\neq 0.
 \end{align}
This leads
\begin{equation}\label{x2}
	x(t,Y)=x(t_0,Y_0)+A_2(Y-Y_0)^5+\mathcal{O}(1)( \lvert t-t_0\rvert^2,\lvert Y-Y_0\rvert^6 )
\end{equation}
So \eqref{S2} can be concluded from \eqref{u2} and \eqref{x2}.
\end{proof}
\begin{rema}
	We found that the generic property of the Novikov equation is different from the property of the Camassa-Holm equation. For the Camassa-Holm equation , if $v(t_0,Y_0)=0$ for some point $(t_0,Y_0)\in\Gamma$ then one can see that $v_t\neq0, v_Y\neq 0 $. Utilizing the implicit function theorem, we can obtain the solution is of $C^2$ smooth in \cite{GRCH}. However, for the Novikov equation, the case $v_Y=0$ will happen, which means we lost some part of the information. We can only assure that the solution is at least differentiable in
	the complement of finitely many characteristic curves.   We believe that this difference is caused by the energy concentration. For the Camassa-Holm equation when the characteristic meet tangentially, they will separate immediately. However, for the Novikov equation, when the characteristics tangentially touch each other, they will stay for a period of time.
\end{rema}

\hfill$\Box$

\smallskip
\noindent\textbf{Acknowledgments} This work was
partially supported by the National Key R\&D Program of China ( No. 2021YFA1002100), the National Natural Science Foundation of China (No.12171493 and No.11701586), 
and the Natural Science Foundation of Guangdong province (No. 2021A1515010296 and 2022A1515011798).


\phantomsection
\addcontentsline{toc}{section}{\refname}
\bibliographystyle{abbrv} 
\bibliography{Feneref11}

\begin{thebibliography}{10}

\bibitem{GRNW}
A.~Bressan and G.~Chen.
\newblock Generic regularity of conservative solutions to a nonlinear wave
  equation.
\newblock {\em Ann. Inst. H. Poincar\'{e} C Anal. Non Lin\'{e}aire},
  34(2):335--354, 2017.

\bibitem{GRDW}
A.~Bressan, T.~Huang, and F.~Yu.
\newblock Structurally stable singularities for a nonlinear wave equation.
\newblock {\em Bull. Inst. Math. Acad. Sin. (N.S.)}, 10(4):449--478, 2015.

\bibitem{CSNW}
A.~Bressan and Y.~Zheng.
\newblock Conservative solutions to a nonlinear variational wave equation.
\newblock {\em Comm. Math. Phys.}, 266(2):471--497, 2006.

\bibitem{NOV}
H.~Cai, G.~Chen, R.~M. Chen, and Y.~Shen.
\newblock Lipschitz metric for the {N}ovikov equation.
\newblock {\em Arch. Ration. Mech. Anal.}, 229(3):1091--1137, 2018.

\bibitem{Cai}
H.~Cai, G.~Chen, Y.~Shen, and Z.~Tan.
\newblock Generic regularity and {L}ipschitz metric for the {H}unter-{S}axton
  type equations.
\newblock {\em J. Differential Equations}, 262(2):1023--1063, 2017.

\bibitem{CWN}
G.~Chen, R.~M. Chen, and Y.~Liu.
\newblock Existence and uniqueness of the global conservative weak solutions
  for the integrable {N}ovikov equation.
\newblock {\em Indiana Univ. Math. J.}, 67(6):2393--2433, 2018.

\bibitem{Chen-Liu}
R.~M. Chen, T.~Hu, and Y.~Liu.
\newblock The shallow-water models with cubic nonlinearity.
\newblock {\em J. Math. Fluid Mech.}, 24(2):Paper No. 49, 31, 2022.

\bibitem{GY}
Z.~Guo, X.~Liu, L.~Molinet, and Z.~Yin.
\newblock Ill-posedness of the {C}amassa-{H}olm and related equations in the
  critical space.
\newblock {\em J. Differential Equations}, 266(2-3):1698--1707, 2019.

\bibitem{HH}
A.~A. Himonas and C.~Holliman.
\newblock The {C}auchy problem for the {N}ovikov equation.
\newblock {\em Nonlinearity}, 25(2):449--479, 2012.

\bibitem{LY}
J.~Li and Z.~Yin.
\newblock Remarks on the well-posedness of {C}amassa-{H}olm type equations in
  {B}esov spaces.
\newblock {\em J. Differential Equations}, 261(11):6125--6143, 2016.

\bibitem{GRCH}
M.~Li and Q.~Zhang.
\newblock Generic regularity of conservative solutions to {C}amassa-{H}olm type
  equations.
\newblock {\em SIAM J. Math. Anal.}, 49(4):2920--2949, 2017.

\bibitem{NE}
V.~Novikov.
\newblock Generalizations of the {C}amassa-{H}olm equation.
\newblock {\em J. Phys. A}, 42(34):342002, 14, 2009.

\bibitem{GRCL}
D.~G. Schaeffer.
\newblock A regularity theorem for conservation laws.
\newblock {\em Advances in Math.}, 11:368--386, 1973.

\bibitem{WY}
X.~Wu and Z.~Yin.
\newblock Well-posedness and global existence for the {N}ovikov equation.
\newblock {\em Ann. Sc. Norm. Super. Pisa Cl. Sci. (5)}, 11(3):707--727, 2012.

\bibitem{WY2}
X.~Wu and Z.~Yin.
\newblock A note on the {C}auchy problem of the {N}ovikov equation.
\newblock {\em Appl. Anal.}, 92(6):1116--1137, 2013.

\bibitem{YL}
W.~Yan, Y.~Li, and Y.~Zhang.
\newblock The {C}auchy problem for the integrable {N}ovikov equation.
\newblock {\em J. Differential Equations}, 253(1):298--318, 2012.

\bibitem{YL2}
W.~Yan, Y.~Li, and Y.~Zhang.
\newblock The {C}auchy problem for the {N}ovikov equation.
\newblock {\em NoDEA Nonlinear Differential Equations Appl.}, 20(3):1157--1169,
  2013.

\bibitem{Yang}
S.~Yang.
\newblock Generic regularity of conservative solutions to the rotational
  {C}amassa-{H}olm equation.
\newblock {\em J. Math. Fluid Mech.}, 22(4):Paper No. 49, 11, 2020.

\end{thebibliography}

\end{document}